\UseRawInputEncoding
\documentclass[reqno]{amsart}
\usepackage{amssymb,amsmath,amsthm}

\usepackage{
  tikz,%
  rotfig, %
  rotfig_tikz, %
  pgfplots, %
  xcolor, %
}
\pgfplotsset{compat=1.18}

\definecolor{SPECorange}{rgb}{1.0,.5625,0}
\definecolor{SPECblue}{rgb}{0,0,0.75}
\definecolor{SPECred}{rgb}{0.75,0,0}
\definecolor{SPECgreen}{rgb}{0,0.75,0}
\definecolor{SPECblack}{rgb}{0.75,0.75,0.75}

\usepackage{arydshln}

\usepackage{siunitx}

\usetikzlibrary{calc}
\usepackage{hyperref}
\hypersetup{colorlinks=true,linkcolor=SPECorange,filecolor=SPECgreen,citecolor=SPECblue,pdfpagemode=UseNone}

\usepackage{subfig}
\usepackage{booktabs}
\usepackage{multirow}
\newcommand{\tablebodyfont}{\small}

\newcounter{mymac@matlab}
\newcommand{\matlab}{MATLAB%
  \ifnum\value{mymac@matlab}<1%
  \textsuperscript{\textregistered}%
  \setcounter{mymac@matlab}{1}%
  \fi%
}

\newif\ifexternalize
\externalizefalse

\title[Fast computation of eigenvalues]{Fast computation of eigenvalues of periodic CMV matrices}

\author[T. Mach]{Thomas Mach$^*$}
\thanks{$^*$Institute for Mathematics, University of Potsdam, Germany.
\texttt{thomas.mach@gmail.com}.}

\author[M. Rinelli and R. Vandebril]{Michele Rinelli$^\dagger$, Raf Vandebril$^\dagger$}
\thanks{$^\dagger$Numerical Analysis and Applied Mathematics (NUMA) unit,
Department of Computer Science, KU Leuven, Leuven, Belgium.
\texttt{michele.rinelli@kuleuven.be}, \texttt{raf.vandebril@kuleuven.be}.}

\author[D. S. Watkins]{David S. Watkins$^\ddagger$}
\thanks{$^\ddagger$Department of Mathematics, Washington State University,
Pullman, WA 99164-3113, USA.
\texttt{davidwatkinsmatrix@gmail.com}.}

\date{}

\keywords{unitary matrix, CMV matrix, Floquet matrix, eigenvalue, fast algorithm}

\subjclass[2020]{65F15, 65H17, 15A18, 65H04}

\begin{document}
\begin{abstract}
Periodic CMV matrices are  unitary matrices that can be specified by $O(n)$ data.  Their eigenvalues can be computed by 
standard methods, storing them as conventional matrices (using $O(n^{2})$ data) in $O(n^{3})$ time.  
Here a fast method that computes the eigenvalues in $O(n^{2})$ time (using $O(n)$ data) is presented.
\end{abstract}

\maketitle

\section{Introduction}  
\label{sec:introduction}

The paper \cite{DaEmFi24} by Damanik, Embree, and Fillman contains an interesting unitary eigenvalue problem.
The starting point is an infinite CMV matrix.   If a certain periodicity condition is satisfied, a related finite 
matrix can be extracted.   This \emph{periodic CMV} matrix (or \emph{Floquet} matrix) 
is an $n \times n$ (with even $n$) unitary matrix 
that is determined by $O(n)$ data, and its eigenvalues are of interest \cite[Prop.\ 2.5]{DaEmFi24}.  
(The eigen\emph{vectors} are not needed.)
Certainly the problem can be solved by storing it 
in a conventional matrix format using $O(n^{2})$ data, and using standard software to compute the eigenvalues in 
$O(n^{3})$ time.   However, one would hope to compute eigenvalues in $O(n^{2})$ time, given that the matrix can
be represented by $O(n)$ data.  This is indeed possible, and we will show how to do it.  

An infinite CMV matrix is a unitary matrix that is a product of two very simple block-diagonal matrices
\begin{displaymath}
L = \left[\begin{array}{cccccc}
\ddots & \\
& C_{0} & \\
& & C_{2} & \\
& & & C_{4} & \\
& & & & \ddots \\
& & & & & 
\end{array}\right]
\end{displaymath}
and
\begin{displaymath}
M = \left[\begin{array}{cccccc}
& & & & & \\
 & \ddots & \\
 & & C_{1} & \\
& & & C_{3} & \\
& & & & C_{5} & \\
& & & & & \ddots 
\end{array}\right].
\end{displaymath}
Each of these matrices is infinite in both directions.  Each block $C_{i}$ is a $2 \times 2$ unitary matrix, so it could be,
for example, a rotator or a reflector.   The blocks in $L$ and $M$ are not aligned.   For example, the first row (column) 
of $C_{1}$ is even with the second row (column) of $C_{0}$, and the second row (column) of $C_{1}$ is even with the 
first row (column) of $C_{2}$.   The product $LM$, which is known as a CMV matrix, is thus not block diagonal but 
pentadiagonal with some special structure.   We will not display this because we are going to keep the matrix in its factored
form.   

Now suppose the matrix entries satisfy a periodicity condition $C_{k+n} = C_{k}$, where $n$ is even.
Then all of the information is contained in a finite submatrix built from $C_{1},\ C_{3}, \ \ldots,  C_{n-1}$
and $C_{2},\ C_{4},\ \ldots,\ C_{n}$, for example, and we can build an $n \times n$ periodic CMV matrix,
the eigenvalues of which are of interest \cite[Prop.\ 2.5]{DaEmFi24}.   

\section{Core transformations}
Before we can proceed, we need to introduce some notation and terminology.   We'll keep it brief;
for more details see, for example, \cite{AuMaRoVaWa18}.

A \emph{core transformation}, briefly \emph{core}, is a unitary matrix that equals the identity matrix except in
two consecutive rows/columns.  Core transformations have lots of uses; we will start by looking at the reduction of
an upper Hessenberg matrix 
\begin{displaymath}
H = \parbox{3.1cm}{
\begin{tikzpicture}[scale=1.66,y=-1cm]
\draw (-.15,-.1) -- (-.2,-.1) -- (-.2,1.1) -- (-.15,1.1);
\draw (1.15,-.1) -- (1.2,-.1) -- (1.2,1.1) -- (1.15,1.1);
\foreach \j in {0,...,5}{
   \foreach \i in {\j,...,5}{\node at (\i/5,\j/5)
     [align=center,scale=1.0]{$\times$};}}
\foreach \j in {0,...,4}{\node at (\j/5,\j/5+.2)
     [align=center,scale=1.0]{$\times$};}
\end{tikzpicture}
}
\end{displaymath}
to upper triangular form.   First a core transformation acts on rows $1$ and $2$ of $H$ to transform the $(2,1)$
entry to zero.   This operation is depicted by
\begin{displaymath}
\parbox{3.1cm}{
\begin{tikzpicture}[scale=1.66,y=-1cm]
\tikzrotation{-.4}{0}
\draw (-.15,-.1) -- (-.2,-.1) -- (-.2,1.1) -- (-.15,1.1);
\draw (1.15,-.1) -- (1.2,-.1) -- (1.2,1.1) -- (1.15,1.1);
\foreach \j in {0,...,5}{
   \foreach \i in {\j,...,5}{\node at (\i/5,\j/5)
     [align=center,scale=1.0]{$\times$};}}
\foreach \j in {0,...,4}{\node at (\j/5,\j/5+.2)
     [align=center,scale=1.0]{$\times$};}
\end{tikzpicture}
}  = \  \
\parbox{2.65cm}{
\begin{tikzpicture}[scale=1.66,y=-1cm]
\phantom{\tikzrotation{0}{0}}
\draw (-.15,-.1) -- (-.2,-.1) -- (-.2,1.1) -- (-.15,1.1);
\draw (1.15,-.1) -- (1.2,-.1) -- (1.2,1.1) -- (1.15,1.1);
\foreach \j in {0,...,5}{
   \foreach \i in {\j,...,5}{\node at (\i/5,\j/5)
     [align=center,scale=1.0]{$\times$};}}
\foreach \j in {1,...,4}{\node at (\j/5,\j/5+.2)
     [align=center,scale=1.0]{$\times$};}
\end{tikzpicture}
}.
\end{displaymath}
The core transformation is denoted by a small double arrow 
\smash{$\!\begin{array}{c}
\Rc \\ \rc
\end{array}\!$}
that points at the first two rows of $H$ to indicate that these are the rows that
are affected by the operation.   If we were to write down the standard matrix representation
of the core, it would look like an identity matrix, except that the first two rows and columns house
a $2 \times 2$ unitary matrix, perhaps a Givens rotation.  

The second step of the transformation applies a core transformation acting on rows $2$ and $3$ to
create a zero in the $(3,2)$ position of the matrix. The result is depicted
as
\begin{displaymath}
\parbox{3.4cm}{
\begin{tikzpicture}[scale=1.66,y=-1cm]
\tikzrotation{-.6}{.2}
\tikzrotation{-.4}{0}
\draw (-.15,-.1) -- (-.2,-.1) -- (-.2,1.1) -- (-.15,1.1);
\draw (1.15,-.1) -- (1.2,-.1) -- (1.2,1.1) -- (1.15,1.1);
\foreach \j in {0,...,5}{
   \foreach \i in {\j,...,5}{\node at (\i/5,\j/5)
     [align=center,scale=1.0]{$\times$};}}
\foreach \j in {0,...,4}{\node at (\j/5,\j/5+.2)
     [align=center,scale=1.0]{$\times$};}
\end{tikzpicture}
}  = \  \
\parbox{2.65cm}{
\begin{tikzpicture}[scale=1.66,y=-1cm]
\phantom{\tikzrotation{0}{0}}
\draw (-.15,-.1) -- (-.2,-.1) -- (-.2,1.1) -- (-.15,1.1);
\draw (1.15,-.1) -- (1.2,-.1) -- (1.2,1.1) -- (1.15,1.1);
\foreach \j in {0,...,5}{
   \foreach \i in {\j,...,5}{\node at (\i/5,\j/5)
     [align=center,scale=1.0]{$\times$};}}
\foreach \j in {2,...,4}{\node at (\j/5,\j/5+.2)
     [align=center,scale=1.0]{$\times$};}
\end{tikzpicture}
}.
\end{displaymath}
This second core also looks like an identity matrix, except for a $2 \times 2$ unitary occupying 
rows and columns $2$ and $3$.   
Continuing the process to its conclusion, we produce upper triangular form:
\begin{displaymath}
\parbox{4.2cm}{
\begin{tikzpicture}[scale=1.66,y=-1cm]
\foreach \j in {0,...,4}{\tikzrotation{-.4-\j/5}{\j/5}}
\draw (-.15,-.1) -- (-.2,-.1) -- (-.2,1.1) -- (-.15,1.1);
\draw (1.15,-.1) -- (1.2,-.1) -- (1.2,1.1) -- (1.15,1.1);
\foreach \j in {0,...,5}{
   \foreach \i in {\j,...,5}{\node at (\i/5,\j/5)
     [align=center,scale=1.0]{$\times$};}}
\foreach \j in {0,...,4}{\node at (\j/5,\j/5+.2)
     [align=center,scale=1.0]{$\times$};}
\end{tikzpicture}
} \ = \  \
\parbox{2.65cm}{
\begin{tikzpicture}[scale=1.66,y=-1cm]
\phantom{\tikzrotation{0}{0}}
\draw (-.15,-.1) -- (-.2,-.1) -- (-.2,1.1) -- (-.15,1.1);
\draw (1.15,-.1) -- (1.2,-.1) -- (1.2,1.1) -- (1.15,1.1);
\foreach \j in {0,...,5}{
   \foreach \i in {\j,...,5}{\node at (\i/5,\j/5)
     [align=center,scale=1.0]{$\times$};}}
\end{tikzpicture}
}.
\end{displaymath}
The core transformations are lined up in what is called an \emph{ascending sequence}.
Inverting the core transformations one after the other, we obtain 
\begin{equation*}
H = \parbox{2.65cm}{
\begin{tikzpicture}[scale=1.66,y=-1cm]
\phantom{\tikzrotation{0}{0}}
\draw (-.15,-.1) -- (-.2,-.1) -- (-.2,1.1) -- (-.15,1.1);
\draw (1.15,-.1) -- (1.2,-.1) -- (1.2,1.1) -- (1.15,1.1);
\foreach \j in {0,...,5}{
   \foreach \i in {\j,...,5}{\node at (\i/5,\j/5)
     [align=center,scale=1.0]{$\times$};}}
\foreach \j in {0,...,4}{\node at (\j/5,\j/5+.2)
    [align=center,scale=1.0]{$\times$};}
\end{tikzpicture}
} \ = \
\parbox{4.7cm}{
\begin{tikzpicture}[scale=1.66,y=-1cm]
\phantom{\tikzrotation{0}{0}}
\foreach \j in {0,...,4}{\tikzrotation{\j/5-1.25}{\j/5}}
\draw (-.15,-.1) -- (-.2,-.1) -- (-.2,1.1) -- (-.15,1.1);
\draw (1.15,-.1) -- (1.2,-.1) -- (1.2,1.1) -- (1.15,1.1);
\foreach \j in {0,...,5}{
   \foreach \i in {\j,...,5}{\node at (\i/5,\j/5)
     [align=center,scale=1.0]{$\times$};}}
\end{tikzpicture}
},
\end{equation*}
In which the (inverted) core transformations form a \emph{descending sequence}.   This equation is a $QR$ decomposition
of $H$.   The product of the cores is a unitary upper Hessenberg matrix, but of course we never form this matrix explicitly;
we keep it in factored form for efficiency.  

\subsection*{Operating with core transformations}

Consider an ascending sequence of three core transformations.
\begin{center}
\begin{tikzpicture}[scale=1.66,y=-1cm]
\node at (-1,0) {$C_{3}C_{2}C_{1} =$};
\node at (0,0) {\begin{tikzpicture}[scale=1.66,y=-1cm]
\foreach \j in {0,1,2}{\tikzrotation{-\j/5}{\j/5}}
\end{tikzpicture}};
\node at (0.5,0) {.};
\end{tikzpicture}
\end{center}
These matrices do not commute as, in general,
$C_{3}C_{2} \neq C_{2}C_{3}$ and $C_{2}C_{1} \neq C_{1}C_{2}$.  
The situation changes, however, if $C_2$ is left out.  $C_1$ and $C_3$ do 
commute since the rows/columns they act on do not overlap: 
\begin{equation}\label{eq:noverlap}
\parbox{7.5cm}{
\begin{tikzpicture}[scale=1.66,y=-1cm]
\node at (-1.8,0) {$C_{3}C_{1} = C_{1}C_{3} = $};
\foreach \j in {0,2}{\tikzrotation{-.4-\j/5}{-.3+\j/5}}
\node at (0,0) {$=$};
\foreach \j in {0,2}{\tikzrotation{.4+\j/5}{-.3+\j/5}}
\node at (1.2,0) {$=$};
\foreach \j in {0,2}{\tikzrotation{1.6}{-.3+\j/5}}
\end{tikzpicture}
}.
\end{equation}
In the end we have placed the symbols for $C_1$ and $C_3$ one atop the other
since the order does not matter.  

Now we consider how to deal with core transformations that act on overlapping 
rows and therefore interact in a nontrivial way.  
In the simplest case we have two
adjacent core transformations acting on the same rows:
\begin{center}
\begin{tikzpicture}[scale=1.66,y=-1cm]
\node at (-.5,0) {$A_{i}B_{i} = $};
\foreach \j in {0,1}{\tikzrotation{\j/5}{-.12}}
\node at (.4,0) {.};
\end{tikzpicture}
\end{center}
In this case we can simply multiply the two core transformations together
to form a single core transformation:  $A_i B_i = C_i$. 
We call this operation \emph{fusion}\index{fusion}.  Pictorially
\begin{center}
\begin{tikzpicture}[scale=1.66,y=-1cm]
\foreach \j in {0,1}{\tikzrotation{\j/5}{-.12}}
\node at (.5,0) {$\rightarrow$};
\tikzrotation{.8}{-.12};
\node at (1,0) {.};
\end{tikzpicture}
\end{center}

Fusions are crucial, but 
a much more frequently used operation is the \emph{turnover}.
Suppose we have three adjacent core transformations $A_{i}B_{i+1}C_{i}$
that don't quite line up:
\begin{displaymath}
A_{i}B_{i+1}C_{i}\,=\parbox{1.2cm}{
\begin{tikzpicture}[scale=1.66,y=-1cm]
\foreach \j in {0,1}{\tikzrotation{\j/5}{\j/5}}
\tikzrotation{.4}{0};
\end{tikzpicture}
}.
\end{displaymath}
If we multiply these three matrices together, the product has action only 
in the three rows and columns $i$, $i+1$, and $i+2$.  It is possible to 
re-factor this matrix in the form $\hat{A}_{i+1}\hat{B}_{i}\hat{C}_{i+1}$, 
thereby ``turning over'' the pattern:

\begin{equation}\label{eq:turnover}
A_{i}B_{i+1}C_{i}\,= \,
\parbox{1.1cm}{\begin{tikzpicture}[scale=1.66,y=-1cm]
\foreach \j in {0,1}{\tikzrotation{\j/5-1.2}{\j/5}}
\tikzrotation{-0.8}{0};
\end{tikzpicture}
}
= \,
\parbox{1.5cm}{\begin{tikzpicture}[scale=1.66,y=-1cm]
\draw (-.15,-.1) -- (-.2,-.1) -- (-.2,0.5) -- (-.15,0.5);
\draw (.55,-.1) -- (.6,-.1) -- (0.6,0.5) -- (.55,0.5);
\foreach \j in {0,1,2}{
   \foreach \i in {0,1,2}{\node at (\i/5,\j/5)
     [align=center,scale=1.0]{$\times$};
   } }
\end{tikzpicture}
}
= 
\parbox{1.1cm}{\begin{tikzpicture}[scale=1.66,y=-1cm]
\tikzrotation{1.3}{.2};
\foreach \j in {0,1}{\tikzrotation{\j/5+1.5}{\j/5}}
\end{tikzpicture}
}
=\, \hat{A}_{i+1}\hat{B}_{i}\hat{C}_{i+1}.
\end{equation}

This operation, which can be done in either direction,  is what we call a \emph{turnover}.  
For information about efficient and accurate implementation of the turnover see \cite[\S 1.4]{AuMaRoVaWa18}.

The turnover operation is usually used to pass a core transformation through a descending or ascending sequence.  
For example, suppose we have a product $C_{1}C_{2}C_{3}C_{4}B_{2}$ composed of a 
descending sequence $C_{1}C_{2}C_{3}C_{4}$ with $B_{2}$ to the right, which
can be pictured as 
\begin{displaymath}
\parbox{1.85cm}{
\begin{tikzpicture}[scale=1.66,y=-1cm]
\foreach \j in {0,...,3}{\tikzrotation{\j/5}{\j/5}}
\tikzrotation{.8}{.2};
\end{tikzpicture}
}
.
\end{displaymath}
We want to somehow pass $B_{2}$ through the descending sequence.  First of 
all, $B_{2}$ commutes with $C_{4}$.  $B_2$ does not commute with $C_{3}$, but we 
can do a turnover:  $C_{2}C_{3}B_{2} = \hat{B}_{3}\hat{C}_{2}\hat{C}_{3}$.  
Finally, $\hat{B}_{3}$ commutes with $C_{1}$, so we have 
\begin{displaymath}
\parbox{1.85cm}{
\begin{tikzpicture}[scale=1.66,y=-1cm]
\foreach \j in {0,...,3}{\tikzrotation{\j/5}{\j/5}}
\tikzrotation{.8}{.2};
\end{tikzpicture}
}
=
\parbox{1.85cm}{
\begin{tikzpicture}[scale=1.66,y=-1cm]
\foreach \j in {0,...,2}{\tikzrotation{\j/5}{\j/5}}
\tikzrotation{.6}{.2};
\tikzrotation{.8}{.6};
\end{tikzpicture}
}
=
\parbox{1.85cm}{
\begin{tikzpicture}[scale=1.66,y=-1cm]
\tikzrotation{0}{0};
\tikzrotation{.2}{.4};
\foreach \j in {1,2,3}{\tikzrotation{.2+\j/5}{\j/5}}
\end{tikzpicture}
} 
=
\parbox{1.85cm}{
\begin{tikzpicture}[scale=1.66,y=-1cm]
\tikzrotation{0}{.4};
\foreach \j in {0,1,2,3}{\tikzrotation{.2+\j/5}{\j/5}}
\end{tikzpicture}
}
.
\end{displaymath}
This is the \emph{shift-through}\index{shift through} operation, and we will usually abbreviate 
it by the simple picture
\begin{equation}\label{op:shiftthroughdown}
\parbox{1.8cm}{
\begin{tikzpicture}[scale=1.66,y=-1cm]
\foreach \j in {0,...,3}{\tikzrotation{\j/5}{\j/5}}
\tikzrotation[gray]{.6}{.2};
\tikzrotation{0}{.4};
\turnoverrl{.6}{.2}{0}{.4};
\end{tikzpicture}
}.
\end{equation}
The core transformation 
$\hat{B}_{3}$ that comes out on the left is, of course, different from the 
core that went in on the right.  It's even been moved down
by one position.  Two of the transformations in the descending 
sequence have been altered as well.  

Clearly we can also pass a core transformation through an ascending sequence:
\begin{displaymath}
 \parbox{1.8cm}{
\begin{tikzpicture}[scale=1.66,y=-1cm]
\foreach \j in {0,...,3}{\tikzrotation{.6-\j/5}{\j/5}}
\tikzrotation{.6}{.4};
\tikzrotation[gray]{0}{.2};
\turnoverlr{0}{.2}{.6}{.4};
\end{tikzpicture}
}.
\end{displaymath}

The one other operation that concerns us is the similarity transformation.   Consider a matrix 
\begin{displaymath}
K = QA = \parbox{3.1cm}{
\begin{tikzpicture}[scale=1.66,y=-1cm]
\tikzrotation{-.4}{.4}
\draw (-.15,-.1) -- (-.2,-.1) -- (-.2,1.1) -- (-.15,1.1);
\draw (1.15,-.1) -- (1.2,-.1) -- (1.2,1.1) -- (1.15,1.1);
\foreach \j in {0,...,5}{
   \foreach \i in {\j,...,5}{\node at (\i/5,\j/5)
     [align=center,scale=1.0]{$\times$};}}
\foreach \j in {0,...,4}{\node at (\j/5,\j/5+.2)
     [align=center,scale=1.0]{$\times$};}
\end{tikzpicture}
}
\end{displaymath}
that is the product of a core transformation $Q$ with some other matrix $A$ (the Hessenberg form depicted here does not matter). 
If we now do a similarity transformation $Q^{-1}KQ$,   we end up with $AQ$.   Thus $Q$ disappears from the left of $A$ and
reappears on the right.  This is depicted 
\begin{displaymath}
\parbox{4.35cm}{
\begin{tikzpicture}[scale=1.66,y=-1cm]
\tikzrotation[gray]{-.4}{.4}
\similaritylr[gray]{-.4}{1.4}{.4}
\tikzrotation{1.4}{.4}
\draw (-.15,-.1) -- (-.2,-.1) -- (-.2,1.1) -- (-.15,1.1);
\draw (1.15,-.1) -- (1.2,-.1) -- (1.2,1.1) -- (1.15,1.1);
\foreach \j in {0,...,5}{
   \foreach \i in {\j,...,5}{\node at (\i/5,\j/5)
     [align=center,scale=1.0]{$\times$};}}
\foreach \j in {0,...,4}{\node at (\j/5,\j/5+.2)
     [align=center,scale=1.0]{$\times$};}
\end{tikzpicture}
}.
\end{displaymath}

\section{Representation of a periodic CMV matrix}

We now return our attention to the infinite CMV matrix $LM$ introduced in Section~\ref{sec:introduction}.
$L$ is a product of core transformations.   This could be represented by a descending sequence, but notice
that these cores all commute with one another, so they can be depicted stacked one atop the other as in 
\ref{eq:noverlap}.   The same is true of the cores in $M$, so 
\begin{displaymath}
LM = \parbox{1.0cm}{
\begin{tikzpicture}[scale=1.2,y=-1cm]
\node at (.1,-.4) [align=center,scale=1.0]{$\scriptstyle{\vdots}$};
\node at (.1,2.2) [align=center,scale=1.0]{$\scriptstyle{\vdots}$};
\foreach \j in {0,...,4}{\tikzrotation{0}{2*\j/5}}
\foreach \j in {0,...,4}{\tikzrotation{.2}{.2+2*\j/5}}
\end{tikzpicture}
}.
\end{displaymath}
When we impose the periodicity condition
$C_{k+n} = C_{k}$, the matrix doubles back on itself.
For illustration we consider the case $n=10$.  We have 
\begin{displaymath}
\parbox{1.0cm}{
\begin{tikzpicture}[scale=1.2,y=-1cm]
\foreach \j in {0,...,4}{\tikzrotation{0}{2*\j/5}}
\foreach \j in {0,...,3}{\tikzrotation{.2}{.2+2*\j/5}}
\tikzrotation[red]{.2}{1.8};
\tikzrotation[red]{.2}{-.2};
\end{tikzpicture}
},
\end{displaymath}
where the two red cores are the same.  To take care of the redundancy, we replace
the two cores by a single spread-out ``core'' that acts on rows/columns $1$ and $n$.
\begin{displaymath}
\parbox{1.1cm}{
\begin{tikzpicture}[scale=1.2,y=-1cm]
\foreach \j in {0,...,4}{\tikzrotation{0}{2*\j/5}}
\foreach \j in {0,...,3}{\tikzrotation{.2}{.2+2*\j/5}}
\tikzrotationhuge[red]{.4}{0}{1.6} 
\end{tikzpicture} }
= \ \ \ \parbox{3.8cm}{
\begin{tikzpicture}[scale=1.2,y=-1cm]
\foreach \j in {0,...,4}{\tikzrotation{0}{2*\j/5}}
\foreach \j in {0,...,3}{\tikzrotation{.2}{.2+2*\j/5}}
\draw (.45,-.1) -- (.4,-.1) -- (.4,1.9) -- (.45,1.9);
\draw (2.5,-.1) -- (2.55,-.1) -- (2.55,1.9) -- (2.5,1.9);
\node at (.55,0) [align=center,scale=1.0]{$\scriptstyle{\overline{c}}$};
\node at (2.35,0) [align=center,scale=1.0]{$\scriptstyle{-s}$};
\node at (.55,1.8) [align=center,scale=1.0]{$\scriptstyle{s}$};
\node at (2.35,1.8) [align=center,scale=1.0]{$\scriptstyle{\phantom{-}c}$};
\foreach \j in {1,...,8}{\node at (.55+\j/5,\j/5)
     [align=center,scale=1.0]{$\scriptstyle{1}$};}.
\end{tikzpicture} }.
\end{displaymath}

This is the periodic CMV matrix whose eigenvalues are wanted.  

We can transform the spread out ``core" into a product of ordinary core transformations by a simple swapping procedure.   
First rows/columns $n$ and $n-1$ are swapped.   Then $n-1$ and $n-2$ are swapped, and so on.  
Finally rows/columns $3$ and $2$ are 
swapped, which leaves the active part of the matrix in rows/columns $1$ and $2$.  This results in 
\begin{displaymath}
\parbox{5.4cm}{
\begin{tikzpicture}[scale=1.2,y=-1cm]
\foreach \j in {0,...,4}{\tikzrotation{0}{2*\j/5}}
\foreach \j in {0,...,3}{\tikzrotation{.2}{.2+2*\j/5}}
\foreach \j in {0,...,8}{\tikzrotation{2.4+\j/5}{\j/5}}
\foreach \j in {1,...,8}{\tikzrotation{2.4-\j/5}{\j/5}}
\end{tikzpicture}
},
\end{displaymath}
where the core transformation at the apex of the pattern is 
$\left[\begin{smallmatrix} \overline{c} & -s \\ s & \phantom{-}c\end{smallmatrix}\right]$, and the others are swapping matrices 
$\left[\begin{smallmatrix} 0 & 1 \\ 1 & 0 \end{smallmatrix}\right]$.   These are reflectors.   If one prefers (as we do) to work with
rotators, one can use rotators  $\left[\begin{smallmatrix} 0 & -1 \\ 1 & \phantom{-}0 \end{smallmatrix}\right]$ on the ascending sequence and the inverses $\left[\begin{smallmatrix} \phantom{-}0 & 1 \\ -1 & 0 \end{smallmatrix}\right]$ on the descending 
sequence.  

The matrix is now represented as a product of $3n-4$ core transformations.   The salient fact 
is that there are $O(n)$ core transformations (hence $O(n)$ data) and not $O(n^{2})$.   Thus it should be
possible to compute the eigenvalues in $O(n^{2})$ time.  We can accomplish this by reducing the matrix to 
upper Hessenberg form, eliminating $O(n)$ core transformations with $O(n^{2})$ work, then computing the
eigenvalues by the $O(n^{2})$ method of \cite{AuMaVaWa17}.

\section{Reduction to upper Hessenberg form}

We begin by eliminating the top and bottom cores of the zig-zag part:
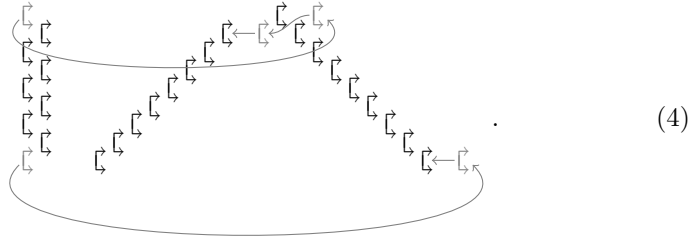
\begin{equation}\label{eq:justthebeginning}
\parbox{7.5cm}{
\begin{tikzpicture}[scale=1.2,y=-1cm]
\foreach \j in {1,...,3}{\tikzrotation{0}{2*\j/5}}
\foreach \j in {0,...,3}{\tikzrotation{.2}{.2+2*\j/5}}
\foreach \j in {0,...,8}{\tikzrotation{2.8+\j/5}{\j/5}}
\foreach \j in {1,...,8}{\tikzrotation{2.4-\j/5}{\j/5}}
\tikzrotation[gray]{0}{0}
\similaritylr[gray]{0}{3.2}{0}
\tikzrotation[gray]{3.2}{0}
\turnoverrl[gray]{3.2}{0}{2.6}{.2}
\tikzrotation[gray]{2.6}{.2}
\shiftthroughrl[gray]{2.6}{2.2}{.2}
\tikzrotation[gray]{0}{1.6}
\similaritylr[gray]{0}{4.8}{1.6}
\tikzrotation[gray]{4.8}{1.6}
\shiftthroughrl[gray]{4.8}{4.4}{1.6}
\end{tikzpicture}
}.
\end{equation}
The bottom core transformation is removed by a similarity transformation followed by a fusion.  The top one is removed 
in the same way, except that a turnover is required before the fusion.

Then the rest of the zigzag part can be pushed into ``the middle" in two steps as follows:
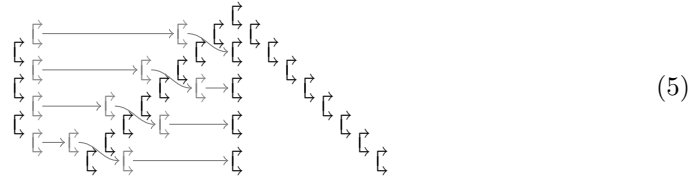
\begin{equation}\label{eq:firstpass}
\parbox{5.5cm}{
\begin{tikzpicture}[scale=1.2,y=-1cm]
\foreach \j in {1,...,3}{\tikzrotation{0}{2*\j/5}}
\foreach \j in {0,...,8}{\tikzrotation{2.4+\j/5}{\j/5}}
\foreach \j in {1,...,8}{\tikzrotation{2.4-\j/5}{\j/5}}
\foreach \j in {0,...,3}{\tikzrotation[gray]{.2}{.2+2*\j/5}}
\shiftthroughlr[gray]{.2}{1.8}{.2}
\tikzrotation[gray]{1.8}{.2}
\turnoverlr[gray]{1.8}{.2}{2.4}{.4}
\shiftthroughlr[gray]{.2}{1.4}{.6}
\tikzrotation[gray]{1.4}{.6}
\turnoverlr[gray]{1.4}{.6}{2.0}{.8}
\tikzrotation[gray]{2.0}{.8}
\shiftthroughlr[gray]{2.0}{2.4}{.8}
\shiftthroughlr[gray]{.2}{1.0}{1.0}
\tikzrotation[gray]{1.0}{1.0}
\turnoverlr[gray]{1.0}{1.0}{1.6}{1.2}
\tikzrotation[gray]{1.6}{1.2}
\shiftthroughlr[gray]{1.6}{2.4}{1.2}
\shiftthroughlr[gray]{.2}{.6}{1.4}
\tikzrotation[gray]{0.6}{1.4}
\turnoverlr[gray]{0.6}{1.4}{1.2}{1.6}
\tikzrotation[gray]{1.2}{1.6}
\shiftthroughlr[gray]{1.2}{2.4}{1.6}

\foreach \j in {0,...,3}{\tikzrotation{2.4}{.4+2*\j/5}}
\end{tikzpicture}
}
\end{equation}

\begin{equation}\label{eq:secondpass}
\parbox{5.5cm}
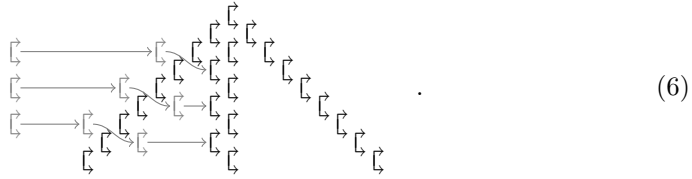
{
\begin{tikzpicture}[scale=1.2,y=-1cm]
\foreach \j in {0,...,8}{\tikzrotation{2.4+\j/5}{\j/5}}
\foreach \j in {1,...,8}{\tikzrotation{2.4-\j/5}{\j/5}}
\foreach \j in {1,...,3}{\tikzrotation[gray]{0}{2*\j/5}}
\shiftthroughlr[gray]{0}{1.6}{.4}
\tikzrotation[gray]{1.6}{.4}
\turnoverlr[gray]{1.6}{.4}{2.2}{.6}
\shiftthroughlr[gray]{0}{1.2}{.8}
\tikzrotation[gray]{1.2}{.8}
\turnoverlr[gray]{1.2}{.8}{1.8}{1.0}
\tikzrotation[gray]{1.8}{1.0}
\shiftthroughlr[gray]{1.8}{2.2}{1.0}
\shiftthroughlr[gray]{0}{0.8}{1.2}
\tikzrotation[gray]{0.8}{1.2}
\turnoverlr[gray]{0.8}{1.2}{1.4}{1.4}
\tikzrotation[gray]{1.4}{1.4}
\shiftthroughlr[gray]{1.4}{2.2}{1.4}

\foreach \j in {1,...,3}{\tikzrotation{2.2}{.2+2*\j/5}}
\foreach \j in {0,...,3}{\tikzrotation{2.4}{.4+2*\j/5}}
\end{tikzpicture}
}.
\end{equation}
So far the cost is $n-2$ turnovers and two fusions, but it turns out that most of the turnovers are free.  Most of the cores in the
ascending sequence are exchange matrices.  The reader can check that in a configuration of the form $C_{i}X_{i+1}X_{i}$, where $X_{i+1}$ and $X_{i}$ are exchange matrices, the result of a turnover is $X_{i+1}X_{i}C_{i+1}$, where $C_{i+1}$ is the same as 
$C_{i}$ but moved down one position.  Thus the turnover costs nothing.   
Initially all of the cores in the ascending sequence except the one at the apex are exchange matrices.   After the operations
of (\ref{eq:justthebeginning}), the top core before the apex ceases to be an exchange matrix.   Thus, in the operations
of (\ref{eq:firstpass}), the top turnover is not free but all of the others are free.   In the process, one more core ceases to 
be an exchange matrix.   In the operations of (\ref{eq:secondpass}), again the top turnover is not free, but all of the others
are.   Thus, of the $n-2$ turnovers so far, all but three are free. 

Now the form of the matrix is 
\begin{equation}\label{eq:tricore_sparse}
\parbox{4.5cm}{
\begin{tikzpicture}[scale=1.2,y=-1cm]
\foreach \j in {0,...,8}{\tikzrotation{2.4+\j/5}{\j/5}}
\foreach \j in {1,...,8}{\tikzrotation{2.4-\j/5}{\j/5}}
\foreach \j in {1,...,3}{\tikzrotation{2.2}{.2+2*\j/5}}
\foreach \j in {0,...,3}{\tikzrotation{2.4}{.4+2*\j/5}}
\end{tikzpicture}
}.
\end{equation}
In \cite[\S~6.2]{AuMaRoVaWa18}  we observed that every unitary matrix can be written as a product of 
$n(n-1)/2$ core transformations in a triangular pattern
\begin{equation}\label{eq:tricore}
\parbox{4.5cm}{
\begin{tikzpicture}[scale=1.2,y=-1cm]
\phantom{\tikzrotation{0}{.4}}
\foreach \i in {0,...,8}{
\foreach \j in {\i,...,8}{\tikzrotation{1.5+\j/5-.4*\i}{\j/5}}}
\end{tikzpicture}
}.
\end{equation}
We also showed how to reduce the matrix to upper Hessenberg form by 
chasing away most of the core transformations, shown in red:
\begin{displaymath}
\parbox{4.5cm}{
\begin{tikzpicture}[scale=1.2,y=-1cm]
\phantom{\tikzrotation{0}{.4}}
\foreach \i in {0,...,0}{
\foreach \j in {\i,...,8}{\tikzrotation[black]{1.5+\j/5-.4*\i}{\j/5}}}
\foreach \i in {1,...,8}{
\foreach \j in {\i,...,8}{\tikzrotation[red]{1.5+\j/5-.4*\i}{\j/5}}}
\end{tikzpicture}
},
\end{displaymath}
leaving only the single descending sequence shown in black.  This is upper Hessenberg form.  
The elimination of the $O(n^{2})$ red cores requires $O(n^{3})$ turnovers, and therefore $O(n^{3})$ work.

The configuration (\ref{eq:tricore_sparse}) is a special case of (\ref{eq:tricore}), so it can be reduced to upper Hessenberg form by the same procedure.   However,  (\ref{eq:tricore_sparse}) is comprised of only $(O(n))$ nontrivial core transformations and is in that sense much sparser than (\ref{eq:tricore}).  
Since the reduction of (\ref{eq:tricore_sparse}) requires the elimination of only $O(n)$ cores, the cost will be $O(n^{2})$. 
The reader can  consult \cite[\S~6.2]{AuMaRoVaWa18} for a general procedure.   
We will demonstrate the reduction in the special case (\ref{eq:tricore_sparse}).   This reduction is different from the 
one shown in \cite{AuMaRoVaWa18}, as it is customized to the case that we are now considering.  
Before we begin, we spread out the picture (\ref{eq:tricore_sparse}) for clarity. 

\begin{equation}\label{cc0} 
\parbox{11cm}{
\begin{tikzpicture}[scale=1.2,y=-1cm]
\phantom{\similaritylr[gray]{.2}{7}{0.6}}
\foreach \j in {0,...,8}{\tikzrotation{5+\j/5}{\j/5}} 
\tikzrotation{4.4}{.2} 
\tikzrotation{4.6}{.4}
\tikzrotation{3.8}{.4} 
\tikzrotation{4.0}{.6}
\tikzrotation{4.2}{.8}
\tikzrotation{3.2}{.6} 
\tikzrotation{2.6}{.8} 
\tikzrotation{2.8}{1.0}
\tikzrotation{3.0}{1.2}
\tikzrotation{2.0}{1.0} 
\tikzrotation{1.4}{1.2} 
\tikzrotation{1.6}{1.4}
\tikzrotation{1.8}{1.6}
\tikzrotation{.8}{1.4} 
\tikzrotation{.2}{1.6}  
\end{tikzpicture}
}
\end{equation}
We will eliminate the extraneous cores from left to right.
Each step begins with a similarity transformation.  The leftmost core is eliminated by a similarity followed by a fusion.
The next one requires a similarity and a turnover, followed by a fusion.   
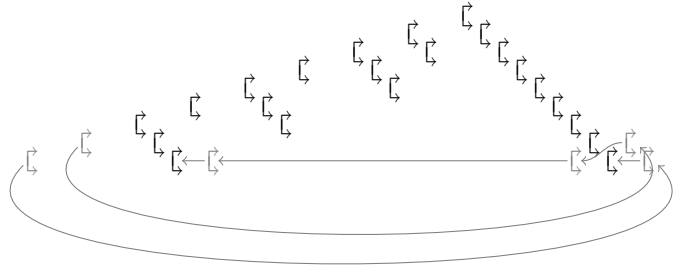
\begin{displaymath} 
\parbox{12cm}{
\begin{tikzpicture}[scale=1.2,y=-1cm]
\foreach \j in {0,...,8}{\tikzrotation{5+\j/5}{\j/5}} 
\tikzrotation{4.4}{.2} 
\tikzrotation{4.6}{.4}
\tikzrotation{3.8}{.4} 
\tikzrotation{4.0}{.6}
\tikzrotation{4.2}{.8}
\tikzrotation{3.2}{.6} 
\tikzrotation{2.6}{.8} 
\tikzrotation{2.8}{1.0}
\tikzrotation{3.0}{1.2}
\tikzrotation{2.0}{1.0} 
\tikzrotation{1.4}{1.2} 
\tikzrotation{1.6}{1.4}
\tikzrotation{1.8}{1.6}
\tikzrotation[gray]{.8}{1.4} 
\similaritylr[gray]{.8}{6.8}{1.4}
\tikzrotation[gray]{6.8}{1.4}
\turnoverrl[gray]{6.8}{1.4}{6.2}{1.6}
\tikzrotation[gray]{6.2}{1.6}
\shiftthroughrl[gray]{6.2}{2.2}{1.6}
\tikzrotation[gray]{2.2}{1.6}
\shiftthroughrl[gray]{2.2}{1.8}{1.6}
\tikzrotation[gray]{.2}{1.6}  
\similaritylr[gray]{.2}{7}{1.6}
\tikzrotation[gray]{7}{1.6}
\shiftthroughrl[gray]{7}{6.6}{1.6}
\end{tikzpicture}
}
\end{displaymath}
Consider next the elimination of the descending sequence of three cores.  
For the first, a turnover is possible after the similarity, but then it gets ``stuck'' 
(red). It will have to wait for elimination on a subsequent step.  The next one follows a similar path and also gets stuck (green).   The one on the bottom is eliminated by fusion.
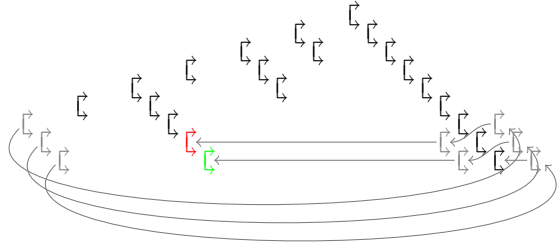
\begin{displaymath} 
\parbox{12cm}{
\begin{tikzpicture}[scale=1.2,y=-1cm]
\phantom{\similaritylr[gray]{.3}{7}{1.6}}
\foreach \j in {0,...,8}{\tikzrotation{5+\j/5}{\j/5}} 
\tikzrotation{4.4}{.2} 
\tikzrotation{4.6}{.4}
\tikzrotation{3.8}{.4} 
\tikzrotation{4.0}{.6}
\tikzrotation{4.2}{.8}
\tikzrotation{3.2}{.6} 
\tikzrotation{2.6}{.8} 
\tikzrotation{2.8}{1.0}
\tikzrotation{3.0}{1.2}
\tikzrotation{2.0}{1.0} 
\tikzrotation[gray]{1.4}{1.2} 
\similaritylr[gray]{1.4}{6.6}{1.2}
\tikzrotation[gray]{6.6}{1.2}
\turnoverrl[gray]{6.6}{1.2}{6.0}{1.4}
\tikzrotation[gray]{6.0}{1.4}
\shiftthroughrl[gray]{6.0}{3.2}{1.4}
\tikzrotation[red]{3.2}{1.4}
\tikzrotation[gray]{1.6}{1.4}
\similaritylr[gray]{1.6}{6.8}{1.4}
\tikzrotation[gray]{6.8}{1.4}
\turnoverrl[gray]{6.8}{1.4}{6.2}{1.6}
\tikzrotation[gray]{6.2}{1.6}
\shiftthroughrl[gray]{6.2}{3.4}{1.6}
\tikzrotation[green]{3.4}{1.6}
\tikzrotation[gray]{1.8}{1.6}
\similaritylr[gray]{1.8}{7}{1.6}
\tikzrotation[gray]{7}{1.6}
\shiftthroughrl[gray]{7}{6.6}{1.6}
\end{tikzpicture}
}
\end{displaymath}
The next solitary core is eliminated by a lengthy path requiring two similarities and three turnovers.
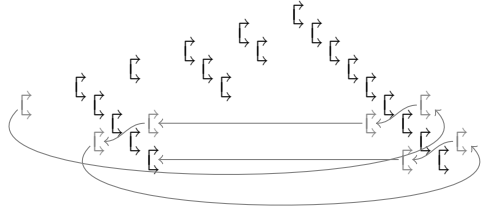
\begin{displaymath} 
\parbox{12cm}{
\begin{tikzpicture}[scale=1.2,y=-1cm]
\phantom{\similaritylr[gray]{.4}{7}{1.6}}
\foreach \j in {0,...,8}{\tikzrotation{5+\j/5}{\j/5}} 
\tikzrotation{4.4}{.2} 
\tikzrotation{4.6}{.4}
\tikzrotation{3.8}{.4} 
\tikzrotation{4.0}{.6}
\tikzrotation{4.2}{.8}
\tikzrotation{3.2}{.6} 
\tikzrotation{2.6}{.8} 
\tikzrotation{2.8}{1.0}
\tikzrotation{3.0}{1.2}
\tikzrotation[black]{3.2}{1.4}
\tikzrotation[black]{3.4}{1.6}
\tikzrotation[gray]{2.0}{1.0} 
\similaritylr[gray]{2.0}{6.4}{1.0}
\tikzrotation[gray]{6.4}{1.0}
\turnoverrl[gray]{6.4}{1.0}{5.8}{1.2}
\tikzrotation[gray]{5.8}{1.2}
\shiftthroughrl[gray]{5.8}{3.4}{1.2}
\tikzrotation[gray]{3.4}{1.2}
\turnoverrl[gray]{3.4}{1.2}{2.8}{1.4}
\tikzrotation[gray]{2.8}{1.4}
\similaritylr[gray]{2.8}{6.8}{1.4}
\tikzrotation[gray]{6.8}{1.4}
\turnoverrl[gray]{6.8}{1.4}{6.2}{1.6}
\tikzrotation[gray]{6.2}{1.6}
\shiftthroughrl[gray]{6.2}{3.4}{1.6}
\end{tikzpicture}
}
\end{displaymath}
The next core, the top one in a descending sequence of five, gets stuck after one turnover (red).  
The next three follow similar paths and also get stuck (green).   We have left the paths out of  the picture to avoid clutter.   
The last core in the group is eliminated in the usual way.
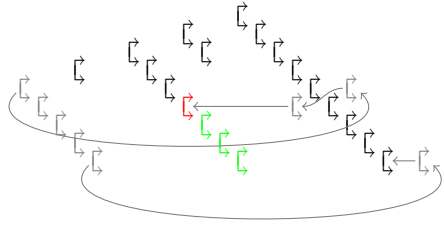
\begin{displaymath} 
\parbox{12cm}{
\begin{tikzpicture}[scale=1.2,y=-1cm]
\phantom{\similaritylr[gray]{.5}{7}{1.6}}
\foreach \j in {0,...,8}{\tikzrotation{5+\j/5}{\j/5}} 
\tikzrotation{4.4}{.2} 
\tikzrotation{4.6}{.4}
\tikzrotation{3.8}{.4} 
\tikzrotation{4.0}{.6}
\tikzrotation{4.2}{.8}
\tikzrotation{3.2}{.6} 
\tikzrotation[gray]{2.6}{.8} 
\similaritylr[gray]{2.6}{6.2}{.8}
\tikzrotation[gray]{6.2}{.8}
\turnoverrl[gray]{6.2}{.8}{5.6}{1.0}
\tikzrotation[gray]{5.6}{1.0}
\shiftthroughrl[gray]{5.6}{4.4}{1.0}
\tikzrotation[red]{4.4}{1.0}
\tikzrotation[gray]{2.8}{1.0}
\tikzrotation[gray]{3.0}{1.2}
\tikzrotation[gray]{3.2}{1.4}
\foreach \j in {1,...,3}{\tikzrotation[green]{4.4+\j/5}{1.0+\j/5}}
\tikzrotation[gray]{3.4}{1.6}
\similaritylr[gray]{3.4}{7}{1.6}
\tikzrotation[gray]{7}{1.6}
\shiftthroughrl[gray]{7}{6.6}{1.6}
\end{tikzpicture}
}
\end{displaymath}
The next solitary core is eliminated by a long path including five turnovers.
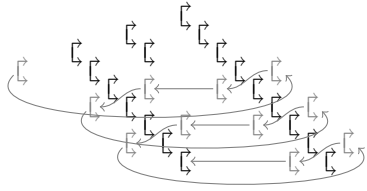
\begin{displaymath} 
\parbox{12cm}{
\begin{tikzpicture}[scale=1.2,y=-1cm]
\phantom{\similaritylr[gray]{.6}{7}{1.6}}
\foreach \j in {0,...,8}{\tikzrotation{5+\j/5}{\j/5}} 
\tikzrotation{4.4}{.2} 
\tikzrotation{4.6}{.4}
\tikzrotation{3.8}{.4} 
\tikzrotation{4.0}{.6}
\tikzrotation{4.2}{.8}
\tikzrotation[black]{4.4}{1.0}
\foreach \j in {1,...,3}{\tikzrotation[black]{4.4+\j/5}{1.0+\j/5}}
\tikzrotation[gray]{3.2}{.6} 
\similaritylr[gray]{3.2}{6.0}{.6}
\tikzrotation[gray]{6.0}{.6}
\turnoverrl[gray]{6.0}{.6}{5.4}{.8}
\tikzrotation[gray]{5.4}{.8}
\shiftthroughrl[gray]{5.4}{4.6}{.8}
\tikzrotation[gray]{4.6}{.8}
\turnoverrl[gray]{4.6}{.8}{4.0}{1.0}
\tikzrotation[gray]{4.0}{1.0}
\similaritylr[gray]{4.0}{6.4}{1.0}
\tikzrotation[gray]{6.4}{1.0}
\turnoverrl[gray]{6.4}{1.0}{5.8}{1.2}
\tikzrotation[gray]{5.8}{1.2}
\shiftthroughrl[gray]{5.8}{5.0}{1.2}
\tikzrotation[gray]{5.0}{1.2}
\turnoverrl[gray]{5.0}{1.2}{4.4}{1.4}
\tikzrotation[gray]{4.4}{1.4}
\similaritylr[gray]{4.4}{6.8}{1.4}
\tikzrotation[gray]{6.8}{1.4}
\turnoverrl[gray]{6.8}{1.4}{6.2}{1.6}
\tikzrotation[gray]{6.2}{1.6}
\shiftthroughrl[gray]{6.2}{5.0}{1.6}
\end{tikzpicture}
}
\end{displaymath}
In the next descending sequence of seven, six get stuck after one turnover, and the last one is eliminated.
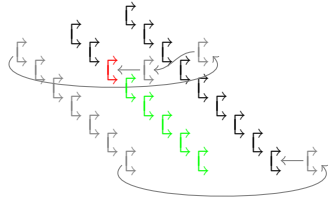
\begin{displaymath} 
\parbox{12cm}{
\begin{tikzpicture}[scale=1.2,y=-1cm]
\phantom{\similaritylr[gray]{.7}{7}{1.6}}
\foreach \j in {0,...,8}{\tikzrotation{5+\j/5}{\j/5}} 
\tikzrotation{4.4}{.2} 
\tikzrotation{4.6}{.4}
\tikzrotation[gray]{3.8}{.4} 
\similaritylr[gray]{3.8}{5.8}{.4}
\tikzrotation[gray]{5.8}{.4}
\turnoverrl[gray]{5.8}{.4}{5.2}{.6}
\tikzrotation[gray]{5.2}{.6}
\shiftthroughrl[gray]{5.2}{4.8}{.6}
\tikzrotation[red]{4.8}{.6}
\foreach \j in {1,...,5}{\tikzrotation[gray]{3.8+\j/5}{.4+\j/5}}
\foreach \j in {1,...,5}{\tikzrotation[green]{4.8+\j/5}{.6+\j/5}}
\tikzrotation[gray]{5.0}{1.6}
\similaritylr[gray]{5.0}{7}{1.6}
\tikzrotation[gray]{7}{1.6}
\shiftthroughrl[gray]{7}{6.6}{1.6}
\end{tikzpicture}
}
\end{displaymath}
Now there is one more descending sequence to be eliminated, and there is no more possiblity 
of getting stuck.   The top core is eliminated by a long path (the longest path!) consisting of four similarity
transforms and seven turnovers.
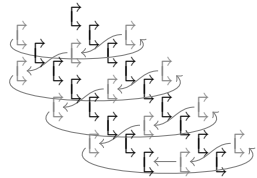
\begin{equation}\label{cc7} 
\parbox{11cm}{
\begin{tikzpicture}[scale=1.2,y=-1cm]
\phantom{\similaritylr[gray]{1.2}{7}{0}}
\foreach \j in {0,...,8}{\tikzrotation{5+\j/5}{\j/5}} 
\tikzrotation[gray]{4.4}{.2} 
\similaritylr[gray]{4.4}{5.6}{.2}
\tikzrotation[gray]{5.6}{.2}
\turnoverrl[gray]{5.6}{.2}{5.0}{.4}
\tikzrotation[gray]{5.0}{.4}
\turnoverrl[gray]{5.0}{.4}{4.4}{.6}
\tikzrotation[gray]{4.4}{.6}
\similaritylr[gray]{4.4}{6.0}{.6}
\tikzrotation[gray]{6.0}{.6}
\turnoverrl[gray]{6.0}{.6}{5.4}{.8}
\tikzrotation[gray]{5.4}{.8}
\turnoverrl[gray]{5.4}{.8}{4.8}{1.0}
\tikzrotation[gray]{4.8}{1.0}
\similaritylr[gray]{4.8}{6.4}{1.0}
\tikzrotation[gray]{6.4}{1.0}
\turnoverrl[gray]{6.4}{1.0}{5.8}{1.2}
\tikzrotation[gray]{5.8}{1.2}
\turnoverrl[gray]{5.8}{1.2}{5.2}{1.4}
\tikzrotation[gray]{5.2}{1.4}
\similaritylr[gray]{5.2}{6.8}{1.4}
\tikzrotation[gray]{6.8}{1.4}
\turnoverrl[gray]{6.8}{1.4}{6.2}{1.6}
\tikzrotation[gray]{6.2}{1.6}
\shiftthroughrl[gray]{6.2}{5.8}{1.6}
\foreach \j in {-1,...,5}{\tikzrotation{4.8+\j/5}{.6+\j/5}}
\end{tikzpicture}
}
\end{equation}
The next one is eliminated by a slightly shorter path 
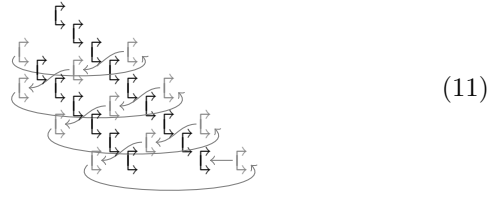
\begin{equation}\label{cc8} 
\parbox{11cm}{
\begin{tikzpicture}[scale=1.2,y=-1cm]
\phantom{\similaritylr[gray]{1}{7}{0}}
\foreach \j in {0,...,8}{\tikzrotation{5+\j/5}{\j/5}} 
\tikzrotation[gray]{4.6}{.4}
\similaritylr[gray]{4.6}{5.8}{.4}
\tikzrotation[gray]{5.8}{.4}
\turnoverrl[gray]{5.8}{.4}{5.2}{.6}
\tikzrotation[gray]{5.2}{.6}
\turnoverrl[gray]{5.2}{.6}{4.6}{.8}
\tikzrotation[gray]{4.6}{.8}
\similaritylr[gray]{4.6}{6.2}{.8}
\tikzrotation[gray]{6.2}{.8}
\turnoverrl[gray]{6.2}{.8}{5.6}{1.0}
\tikzrotation[gray]{5.6}{1.0}
\turnoverrl[gray]{5.6}{1.0}{5.0}{1.2}
\tikzrotation[gray]{5.0}{1.2}
\similaritylr[gray]{5.0}{6.6}{1.2}
\tikzrotation[gray]{6.6}{1.2}
\turnoverrl[gray]{6.6}{1.2}{6.0}{1.4}
\tikzrotation[gray]{6.0}{1.4}
\turnoverrl[gray]{6.0}{1.4}{5.4}{1.6}
\tikzrotation[gray]{5.4}{1.6}
\similaritylr[gray]{5.4}{7.0}{1.6}
\tikzrotation[gray]{7.0}{1.6}
\shiftthroughrl[gray]{7.0}{6.6}{1.6}
\foreach \j in {0,...,5}{\tikzrotation{4.8+\j/5}{.6+\j/5}}
\end{tikzpicture}
}
\end{equation}
The eliminations continue in this way, until the final elimination, which requires only a similarity and a fusion, after which the
matrix is upper Hessenberg.   
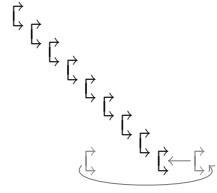
\begin{displaymath} 
\parbox{12cm}{
\begin{tikzpicture}[scale=1.2,y=-1cm]
\phantom{\similaritylr[gray]{1.2}{7}{1.2}}
\foreach \j in {0,...,8}{\tikzrotation{5+\j/5}{\j/5}} 
\tikzrotation[gray]{5.8}{1.6}
\similaritylr[gray]{5.8}{7.0}{1.6}
\tikzrotation[gray]{7.0}{1.6}
\shiftthroughrl[gray]{7.0}{6.6}{1.6}
\end{tikzpicture}
}
\end{displaymath}
Now the eigenvalues of the upper Hessenberg matrix can be computed with $O(n^{2})$ work by the 
method of \cite{AuMaVaWa17}.

We now make some observations about the work required for the reduction procedure.  
Starting from (\ref{cc0}),  $2n-5$ core transformations are removed.  
In each case the removal is done by fusion, so there are $2n-5$ fusions in all.   These 
take place at the bottom.  The mechanism that moves cores toward the bottom is the turnover.
To move the uppermost core to the bottom requires $n-3$ turnovers, as illustrated in (\ref{cc7}).  All other cores 
require fewer than $n-3$, as shown in (\ref{cc8}), for example.  In each case the number of turnovers is at most
$O(n)$, and since there are $O(n)$ cores, and each turnover requires $O(1)$ work, the total work is $O(n^{2})$.  
The reader can check that the exact number of turnovers, starting from (\ref{cc0}), is $(n-3)^{2}$.

\section{Numerical Experiments}

We implemented the reduction to Hessenberg form in Fortran. The method of \cite{AuMaVaWa17} computes the eigenvalues of the resulting unitary Hessenberg matrix; the corresponding Fortran implementation is included in the EISCOR package, available from~\url{https://github.com/eiscor/}. Our implementation of the method described in this paper, together with the code used for all numerical experiments reported below, is available at~\url{https://numa.cs.kuleuven.be/software/momentum-software}.

We found that the reduction phase takes about half as much time as the eigenvalue computation phase.  

We compared our code against standard codes that ignore the unitary structure.   These included
the LAPACK Fortran routines ZGEEV and ZLAHQR, and MATLAB's \texttt{eig} function.  
In our computing environment \texttt{eig} is fastest, so that's what we are showing here.  
Randomly generated unitary CMV matrices were used in these experiments.  

The experiments were run on a computer with an AMD Ryzen 7 PRO 7840U processor
and 16 GB of memory.  The Fortran codes were compiled with GNU Fortran 15.2.0
using the optimization flag \texttt{-O3}. 

We display our results in both graphic and tabular form.  
The graph shows computing time as a function of matrix dimension for $n=10$ to 5000.   We
see that our code (blue) is faster than the structure-ignoring code (red) 
for all values of $n$.   At $n=10$ the improvement is modest, but as $n$ grows the improvement becomes
substantial.  

The same results are shown in tabular form, which also lists the speedup and maximum errors.   
At $n = 100$ the speedup is already about 5, and at $n=5000$ it is over 30.  

These methods are both backward stable, and the eigenvalues of unitary 
matrices are well conditioned, so the results are guaranteed to be accurate.   
Nevertheless we checked the our results against ``exact'' values computed by the ADVANPIX Multiprecision 
Computing Toolbox~\cite{advanpix2015} in quadruple precision.   For each value of $n$ we
did ten runs with different randomly generated matrices; the maximum errors reported in the table are maxima
over all ten runs.   We see that all of the errors are tiny, as expected.   The errors for our structured code are 
slightly smaller than for the unstructured code.  


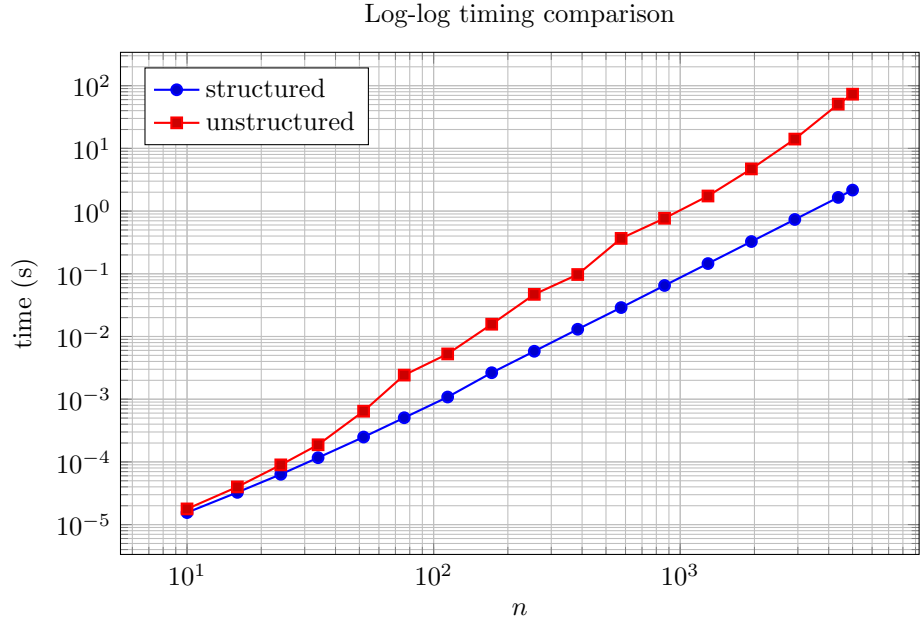
\begin{center}
\begin{tikzpicture}
\begin{axis}[
  width=0.96\textwidth,
  height=0.65\textwidth,
  xmode=log,
  ymode=log,
  log basis x=10,
  log basis y=10,
  grid=both,
  xlabel={$n$},
  ylabel={time (s)},
  legend cell align=left,
  legend pos=north west,
  title={Log-log timing comparison}
]
\addplot+[mark=*, thick, blue]
  table[col sep=comma, x=n, y=feiscor_total]
  {random_feiscor_vs_matlab_eig.csv};
\addlegendentry{structured}

\addplot+[mark=square*, thick, red]
  table[col sep=comma, x=n, y=matlab_eig]
  {random_feiscor_vs_matlab_eig.csv};
\addlegendentry{unstructured}
\end{axis}
\end{tikzpicture}
\end{center}

\par\addvspace{8pt}

\begin{center}
\tablebodyfont
\begin{tabular}{rllcrr}
\toprule
& \multicolumn{3}{c}{time (seconds)}  & \multicolumn{2}{c}{maximum error} \\
\midrule
$n$ & structured &  unstructured & speedup &  structured & unstructured\\
\midrule
10 & \num{1.554e-05} & \phantom{0}\num{1.780e-05} & \phantom{0}1.1x & \num{1.4e-15} & \num{1.7e-15} \\ 
16 & \num{3.265e-05} & \phantom{0}\num{3.994e-05} & \phantom{0}1.2x & \num{1.5e-15}  & \num{3.6e-15} \\
24 & \num{6.342e-05} & \phantom{0}\num{8.968e-05} & \phantom{0}1.4x & \num{1.7e-15}  & \num{3.6e-15} \\
34 &  \num{1.167e-04} & \phantom{0}\num{1.875e-04} & \phantom{0}1.6x & \num{1.7e-15}  & \num{3.6e-15} \\
52 & \num{2.492e-04} & \phantom{0}\num{6.437e-04} & \phantom{0}2.6x & \num{3.4e-15}  & \num{6.1e-15} \\
76 & \num{5.063e-04} & \phantom{0}\num{2.424e-03} & \phantom{0}4.8x & \num{2.9e-15}  & \num{8.7e-15} \\
114 & \num{1.083e-03} & \phantom{0}\num{5.273e-03} & \phantom{0}4.9x & \num{4.2e-15}  & \num{9.9e-15} \\
172 & \num{2.644e-03} & \phantom{0}\num{1.570e-02} & \phantom{0}5.9x & \num{5.1e-15}  & \num{1.3e-14} \\
256 & \num{5.812e-03} & \phantom{0}\num{4.697e-02} & \phantom{0}8.1x & \num{5.8e-15}  & \num{1.5e-14} \\
384 & \num{1.304e-02} & \phantom{0}\num{9.736e-02} & \phantom{0}7.5x & \num{1.0e-14}  & \num{2.0e-14} \\
576 & \num{2.893e-02} & \phantom{0}\num{3.668e-01} & 12.7x & \num{1.2e-14}  & \num{4.9e-14} \\
864 & \num{6.508e-02} & \phantom{0}\num{7.680e-01} & 11.8x & \num{1.4e-14}  & \num{5.3e-14} \\
1296 & \num{1.455e-01} & \phantom{0}1.741 & 12.0x & \num{2.0e-14}  & \num{4.4e-14} \\
1944 & \num{3.268e-01} & \phantom{0}4.707 & 14.4x & \num{2.7e-14}  & \num{5.7e-14} \\
2916 & \num{7.349e-01} & 14.068 & 19.1x & \num{4.2e-14}  & \num{5.0e-14} \\
4374 & 1.652 & 50.740 & 30.7x & \num{5.3e-14}  & \num{9.1e-14} \\ 
5000 & 2.162 & 73.181 & 33.8x & \num{6.1e-14}  & \num{9.1e-14} \\ 
\bottomrule
\end{tabular}
\end{center}

\section*{Funding}
The research was supported by the Research Foundation -- Flanders (FWO), via the junior postdoctoral fellowship 12A1325N (Short Recurrences for Block Krylov Methods with Applications to Matrix Functions and Model Order Reduction) for the second author, projects
G0A9923N (Low rank tensor approximation techniques for up- and downdating of massive online time series clustering) and G0B0123N (Short recurrence relations for rational Krylov and orthogonal rational functions inspired by modified moments) for the third author. The third author is also supported by the Research Council KU Leuven (Belgium), project C16/21/002 (Manifactor: Factor Analysis for Maps into Manifolds).

\bibliographystyle{amsplain}
\bibliography{periodic_un}
\end{document}